\documentclass[12pt]{amsart}

\usepackage{amscd,amssymb}
\usepackage{amsthm,amsmath,amssymb}

\newtheorem{theorem}[equation]{Theorem}

\newtheorem{proposition}[equation]{Proposition}
\newtheorem{lemma}[equation]{Lemma}
\newtheorem{corollary}[equation]{Corollary}

\theoremstyle{definition}
\newtheorem{definition}[equation]{Definition}

\theoremstyle{remark}
\newtheorem{remark}[equation]{Remark}
\newtheorem{example}[equation]{Example}

\title{Double spaces with isolated singularities}

\author{Ivan Cheltsov}

\address{Steklov Institute, 8 Gubkin street, Moscow 117966, Russia} %

\thanks{The author is very grateful to A.Corti, M.Grinenko,
V.Is\-kov\-skikh, S.Kudryavtsev, V.Kulikov, M.Mel\-la, J.~Park,
Yu.Pro\-kho\-rov, A.Pukh\-li\-kov, V.Sho\-ku\-rov and L.Wotzlaw
for fruitful conversations.}

\email{cheltsov@yahoo.com}

\begin{document}

\begin{abstract}
We prove the non-rationality of a double cover of $\mathbb{P}^{n}$
branched over a hypersurface $F\subset\mathbb{P}^{n}$ of degree
$2n$ having isolated singularities such that $n\ge 4$ and every
singular points of the hypersurface $F$ is ordinary, i.e. the
projectivization of its tangent cone is smooth, whose
multiplicity does not exceed $2(n-2)$.%
\end{abstract}

\maketitle

\section{Introduction.}
\label{section:introducation}

For a given algebraic variety, it is one of the most substantial
questions whether it is rational\footnote{A variety is called
rational when it is birationally isomorphic to $\mathbb{P}^{n}$,
i.e. when its field of rational functions is a purely
transcendental extension of the base field.} or not. Global
holomorphic differential forms are natural birational invariants
of a smooth algebraic variety that solve the rationality problem
for algebraic curves and surfaces (see \cite{Za58}). However,
there are only four known methods to prove the non-rationality of
a rationally connected higher-dimensional (see \cite{Is97}). In
the following we assume that all varieties are projective, normal,
and defined over $\mathbb{C}$.

The non-rationality of a smooth quartic 3-fold was proved in
\cite{IsMa71} using the group of birational automorphisms as a
birational invariant. The non-rationality of a smooth cubic 3-fold
was proved in \cite{ClGr72} through the study of its intermediate
Jacobian. The birational invariance of the torsion subgroup of the
group $H^{3}(\mathbb{Z})$ was used in \cite{ArMu72} to prove the
non-rationality of some unirational varieties. The non-rationality
of a wide class of rationally connected varieties was proved in
\cite{Ko95} using the degeneration technique and the reduction
into positive characteristic (see \cite{Ko96}, \cite{Ko00},
\cite{ChWo04}).

\begin{definition}
\label{definition:super-rigidity} A terminal
$\mathbb{Q}$-factorial Fano variety $V$ with $\mathrm{Pic}(V)\cong
\mathbb{Z}$ is called bi\-ra\-ti\-o\-nal\-ly super-rigid if the
following three conditions hold: the variety $V$ is not birational
to a fibration\footnote{For a fibration $\tau:Y\to Z$ we assume
$\dim(Y)>\dim(Z)\ne 0$ and
$\tau_{*}(\mathcal{O}_{Y})=\mathcal{O}_{Z}$.} whose general fiber
has Kodaira dimension $-\infty$; the variety $V$ is not birational
to a Fano variety with terminal $\mathbb{Q}$-factorial
singularities, whose Picard group is $\mathbb{Z}$ and that is not
biregular to $V$; the groups $\mathrm{Bir}(V)$ and
$\mathrm{Aut}(V)$ coincide.
\end{definition}

The notion of birational super-rigidity goes back to
\cite{IsMa71}. For example, the paper \cite{IsMa71} implicitly
proves that any smooth quartic 3-fold in $\mathbb{P}^4$ is
birationally super-rigid (see \cite{Co95}).

\begin{remark}
\label{remark:rigidity-vs-rationality} A birationally super-rigid
Fano variety is not rational and not birational to a conic bundle.
However, there are non-rational Fano varieties that are not
birationally super-rigid, e.g. a smooth cubic 3-fold.
\end{remark}

Let $\pi:X\to \mathbb{P}^{n}$ be a double cover branched over a
hypersurface $F$ of degree $2n$ with isolated singularities. Then
$K_{X}\sim\pi^{*}(\mathcal{O}_{\mathbb{P}^{n}}(-1))$. So, the
variety $X$ is a Fano variety.

\begin{remark}
\label{remark:Kollar} The variety $X$ is known to be birationally
super-rigid in the following three cases: $n\ge 3$ and $F$ is
smooth (see \cite{Is80}, \cite{Pu88a}); $n\ge 3$ and the
hypersurface $F$ has one ordinary singular point of even
multiplicity that does not exceed $2(n-2)$ (see \cite{Pu97});
$n=3$ and the variety $X$ is nodal and $\mathbb{Q}$-factorial (see
\cite{ChPa04}). For $n\ge 3$ the non-rationality of a double cover
of $\mathbb{P}^{n}$ ramified in a very general hypersurface of
degree greater than ${\frac {n+1} {2}}$ is proved in \cite{Ko96}.
\end{remark}

The main purpose of this paper is to prove the following result.

\begin{theorem}
\label{theorem:main} Suppose that $n\ge 4$ and every singular
point $O$ of $F$ is ordinary, i.e. the projectivization of a
tangent cone to $F$ at $O$ is smooth, such that
$\mathrm{mult}_{O}(F)\le 2(n-2)$. Then $X$ is birationally
super-rigid.
\end{theorem}

\begin{corollary}
\label{corollary:Bir} In the conditions of
Theorem~\ref{theorem:main}, the group $\mathrm{Bir(X)}$ is finite.
\end{corollary}

\begin{corollary}
\label{corollary:nodal} A double cover of $\mathbb{P}^{n}$
branched over a nodal hypersurface of degree $2n$ with any number
of ordinary double points is not birationally equivalent to any
elliptic fibration for $n\ge 4$.
\end{corollary}

\begin{example}
\label{example:nodal-plane-odd} Let $n=2k$ for $k\in\mathbb{N}$
and $F\subset\mathbb{P}^{2k}$ be a sufficiently general
hypersurface of degree $4k$ passing through a linear subspace
$\mathrm{\Pi}\subset\mathbb{P}^{2k}$ of dimension $k$. The variety
$X$ can be given by the equation
$$
y^{2}=\sum_{i=1}^{k}a_{i}(x_{0},\ldots,x_{2k+1})x_{i}\subset\mathbb{P}(1^{2k+1},2k)\cong\mathrm{Proj}(\mathbb{C}[x_{0},\ldots,x_{2k+1},y]),%
$$
where $a_{i}$ is a homogeneous polynomial of degree $4k-1$, and
the linear subspace $\mathrm{\Pi}\subset\mathbb{P}^{n}$ is given
by $x_{1}=\ldots=x_{k}=0$. The hypersurface $F$ is nodal, it has
$(4k-1)^{k}$ ordinary double points given by the equations
$$
a_{1}=\ldots=a_{k}=x_{1}=\ldots=x_{k}=0,
$$
and $X$ is non-rational for $k\ge 2$ by
Corollary~\ref{corollary:nodal}.
\end{example}

\begin{example}
\label{example:nodal-plane-even} Let $n=2k+1$ for $k\in\mathbb{N}$
and $F\subset\mathbb{P}^{2k+1}$ be a sufficiently general
hypersurface of degree $4k+2$ that is given by the equation
$$
g^{2}(x_{0},\ldots,x_{2k+2})=\sum_{i=1}^{k}a_{i}(x_{0},\ldots,x_{2k+2})b_{i}(x_{0},\ldots,x_{2k+2}),%
$$
where $g$, $a_{i}$ and $b_{i}$ are homogeneous polynomials of
degree $2k+1$, and $x_{i}$ is a homogeneous coordinate on
$\mathbb{P}^{n}$. The hypersurface $F$ is nodal, and it has
$(2k+1)^{2k+1}$ ordinary double points given by the equations
$$
g=a_{1}=\ldots=a_{k}=b_{1}=\ldots=b_{k}=0,
$$
and the double cover $\pi:X\to\mathbb{P}^{2k+1}$ branched over $F$
is non-rational and birationally super-rigid for $k\ge 2$ by
Theorem~\ref{theorem:main}. In the case $k=1$ one can unproject
(see \cite{Re00}) the variety $X$ into a fibration of cubic
surfaces, i.e. the variety $X$ is not birationally super-rigid. In
the latter case it is unknown whether $X$ is rational or not.
\end{example}

\begin{remark}
\label{remark:Varchenko} In the conditions of
Theorem~\ref{theorem:main}, the best known upper bound of the
number of ordinary singular points of the hypersurface
$F\subset\mathbb{P}^{n}$ is due to \cite{Va83}. Namely,
$|\mathrm{Sing}(F)|\le\mathrm{A}_{n}(2n)$, where
$\mathrm{A}_{n}(2n)$ is a number of integer points
$(a_{1},\ldots,a_{n})\subset\mathbb{R}^{n}$ such that
$(n-1)^{2}<\sum_{i=1}^{n}a_{i}\le n^{2}$ and all $a_{i}\in(0,2n)$.
Hence, $|\mathrm{Sing}(X)|$ does not exceed $68$, $1190$ and
$27237$ when $n=3$, $4$ and $5$ respectively. It is expected that
this bound is far from being sharp for $n\gg 0$ (cf.
\cite{vSt93}). In the case $n=3$ there is a sharp bound
$|\mathrm{Sing}(X)|\le 65$ (see \cite{St78}, \cite{CaCe82},
\cite{Ba96}, \cite{JaRu97}, \cite{Wa98}).
\end{remark}

The condition $\mathrm{mult}_{O}(F)\le 2(n-2)$ in
Theorem~\ref{theorem:main} can not be omitted.

\begin{example}
\label{example:conic-bundle} Let $O$ be a singular point of  $F$
such that $\mathrm{mult}_{O}(F)=2(n-1)$, and
$\gamma:\mathbb{P}^{n}\dashrightarrow\mathbb{P}^{n-1}$ be a
projection from $O$. Then the normalization of the general fiber
of $\gamma\circ\pi$ is a smooth rational curve, i.e. $X$ is
birationally isomorphic to a conic bundle.
\end{example}

The condition $n\ge 4$ in Theorem~\ref{theorem:main} can not be
omitted.

\begin{example}
\label{example:Barth} Let $n=3$ and $F$ be a Barth sextic (see
\cite{Ba96}) given by
$$
4(\alpha^2x^2-y^2)(\alpha^2y^2-z^2)(\alpha^2z^2-x^2)-t^2(1+2\alpha)(x^2+y^2+z^2-t^{2})^2=0
$$
in $\mathbb{P}^{3}\cong \mathrm{Proj}(\mathbb{C}[x,y,z,t])$, where
$\alpha={\frac {1+\sqrt{5}} {2}}$. Then $X$ has only ordinary
double points, $|\mathrm{Sing}(X)|=65$, and $X$ is birational to a
determinantal quartic 3-fold in $\mathbb{P}^4$ with $42$ nodes
(see \cite{En99}, \cite{Pet98}). Thus, $X$ is rational.
\end{example}

The claim of Theorem~\ref{theorem:main} holds for $n=3$ in the
additional assumption that $X$ is $\mathbb{Q}$-factorial, which is
always the case when the number of nodes of $X$ does not exceed
$14$ due to \cite{ChPa04}. On the other hand, there are nodal
double covers of $\mathbb{P}^{3}$ with $15$ nodes that are not
$\mathbb{Q}$-factorial and not birationally super-rigid

The nature of Theorem~\ref{theorem:main} is a reminiscence of the
Noether theorem on the structure of the group
$\mathrm{Bir}(\mathbb{P}^{2})$ (see \cite{Ne71}, \cite{Is80},
\cite{Co95}). The relevant problem is to classify pencils of plane
elliptic curves up to the action of the group
$\mathrm{Bir}(\mathbb{P}^{2})$. It was studied in \cite{Be77}. The
ideas of \cite{Be77} were recovered later in \cite{Do66}, where it
was proved that any pencil of plane elliptic curves can be
birationally transformed into a special elliptic pencil, so-called
Halphen pencil (see \cite{Hal82} and \S 5.6 of \cite{CosDo89}).
The similar problem can be considered for the variety $X$ as well.
Namely, we prove the following result.

\begin{theorem}
\label{theorem:second} In the conditions of
Theorem~\ref{theorem:main}, let $\rho:X\dashrightarrow Z$ be a
rational map such that the normalization of a general fiber of
$\rho$ is a connected elliptic curve. Then there is a point $O$ of
the hypersurface $F$ and a birational map
$\gamma:\mathbb{P}^{n-1}\dasharrow Y$ such that
$\mathrm{mult}_{O}(F)=2(n-2)$ and $\rho=\gamma\circ\beta\circ\pi$,
where $\beta:\mathbb{P}^{n}\dasharrow\mathbb{P}^{n-1}$ is a
projection from the point $O\in\mathbb{P}^{n}$.%
\end{theorem}

\begin{example}
\label{example:elliptic} Let $F\subset\mathbb{P}^{n}$ be a
hypersurface given by the equation
$$
\sum_{i=0}^{4}g_{2n-i}(x_{1},\ldots,x_{n})x_{0}^{i}=0\subset\mathbb{P}^{n}\cong\mathrm{Proj}(\mathbb{C}[x_{0},\ldots,x_{n}]),%
$$
where $g_{i}$ is a general homogeneous polynomial of degree $i$.
The hypersurface $F$ is smooth outside a point
$(1:0:\cdots:0)\in\mathbb{P}^{n}$, which is an ordinary singular
point of $F$ of multiplicity $2n-2$. Thus, in the case $n\ge 4$
the variety $X$ is birationally equivalent to a single elliptic
fibration induced by the projection from the point $O$ by
Theorem~\ref{theorem:second}.
\end{example}

\begin{corollary}
\label{corollary:elliptic-nodal} A double cover of
$\mathbb{P}^{n}$ branched over a nodal hypersurface of degree $2n$
with any number of ordinary double points is not birationally
equivalent to any elliptic fibration for $n\ge 4$.
\end{corollary}

The condition $n\ge 4$ in Theorem~\ref{theorem:second} can not be
omitted (see \cite{ChPa04}).

\begin{example}
\label{example:special-elliptic-fibration} Let $n=3$ and
$F\subset\mathbb{P}^3$ be a nodal sextic such that $F$ contains a
line $L\subset\mathbb{P}^3$ and the set $\mathrm{Sing}(F)\cap L$
consists of $4$ nodes. For a sufficiently general point $P\in X$,
there is a unique hyperplane $H\subset\mathbb{P}^3$ passing
through the point $\pi(P)$ and the line $L$. For a quintic curve
$C\subset H$ given by $F\cap H=L\cup C$, the intersection
$L\cap(C\setminus\mathrm{Sing}(X))$ consists of a single point
$Q$. Take a line $L_{P}\subset \mathbb{P}^3$ passing through
$\pi(P)$ and $Q$ and define a rational map
$\mathrm{\Xi}:X\dasharrow \mathrm{Gr}(2,4)$ by
$\mathrm{\Xi}(P)=L_{P}$. The normalization of a general fiber of
the map $\mathrm{\Xi}$ is an elliptic curve. The rational map
$\mathrm{\Xi}$ can not be obtained by means of the construction in
Theorem~\ref{theorem:second}.
\end{example}

Birational transformations of smooth Fano 3-folds into elliptic
fibrations were used in \cite{BoTsch99}, \cite{BoTsch00},
\cite{HaTsch00} in the proof of the following result.

\begin{theorem}
\label{theorem:potential-density} The set of rational points is
potentially dense\footnote{The set of rational points of a variety
$V$ defined over a number field $\mathbb{F}$ is called potentially
dense if for a finite extension of fields
$\mathbb{K}\slash\mathbb{F}$ the set of $\mathbb{K}$-rational
points of the variety $V$ is Zariski dense.} on every smooth Fano
3-fold defined over a number field $\mathbb{F}$ with a possible
exception of a double cover of $\mathbb{P}^{3}$ ramified in a
smooth sextic surface.
\end{theorem}

The possible exception appears in
Theorem~\ref{theorem:potential-density} because a smooth sextic
double solid is the only smooth Fano 3-fold that is not
birationally isomorphic to an elliptic fibration (see
\cite{Ch00a}). For results relevant to
Theorem~\ref{theorem:second} see \cite{Ch00a}, \cite{Ch00b},
\cite{Ch00c}, \cite{Ch03b}, \cite{Ch04a}, and \cite{Ry02}. The
proof of Theorem~\ref{theorem:second} implicitly gives the
following result.

\begin{theorem}
\label{theorem:third} In the conditions of
Theorem~\ref{theorem:main}, the variety $X$ is not birationally
isomorphic to a Fano variety with canonical singularities.
\end{theorem}

The condition $n\ge 4$ in Theorem~\ref{theorem:third} can not be
omitted (see \cite{ChPa04}).

\section{Preliminary results.}
\label{section:preliminaries}

In this chapter we consider properties of usual log-pairs and
so-called movable log pairs (see \cite{Al91}, \cite{Ch00a},
\cite{Ch03b}). The basic notions, notations and definitions can be
found in \cite{KMM}, \cite{Ko91}, \cite{Co95}, \cite{Co00},
\cite{Pu00}, \cite{Ch03b}. A priori we do not assume any
restriction on the coefficients of the considered boundaries.

\begin{theorem}
\label{theorem:Nother-Fano} Let $X$ be a Fano variety having
terminal $\mathbb{Q}$-factorial singularities and
$\mathrm{Pic}(X)\cong \mathbb{Z}$ such that the set of centers of
canonical singularities  $\mathbb{CS}(X, M_{X})$ is empty for
every movable log pair $(X, M_{X})$ such that $M_{X}$ is effective
and
 $-(K_{X}+M_{X})$ is ample. Then $X$ is birationally
super-rigid.
\end{theorem}

\begin{proof}
See \cite{Co95}, \cite{Pu00} or \cite{Ch03b}.
\end{proof}

\begin{theorem}
\label{theorem:elliptic-fibrations} Let $X$ be a Fano variety with
terminal $\mathbb{Q}$-factorial singularities and
$\mathrm{Pic}(X)\cong \mathbb{Z}$, $\rho:X\dasharrow Y$ be a
birational map, $\tau:Y\to Z$ be a fibration whose general fiber
has Kodaira dimension zero, $H$ be a very ample divisor on $Z$,
and $M_{X}=r\rho^{-1}(|\tau^{*}(H)|)$ for a positive rational
number $r$ such that $K_{X}+M_{X}\sim_\mathbb{Q} 0$. Then the set
of centers of canonical singularities $\mathbb{CS}(X, M_{X})$ is
not empty.
\end{theorem}

\begin{proof}
See \cite{Ch00a}, \cite{Ch03b} and \cite{ChPa04}.
\end{proof}

\begin{theorem}
\label{theorem:canonical-Fanos} Let $X$ be a Fano variety with
terminal $\mathbb{Q}$-factorial singularities and
$\mathrm{Pic}(X)\cong \mathbb{Z}$, $\rho:X\dasharrow Y$ be a
non-biregular birational map, $Y$ be a Fano variety with canonical
singularities, and $M_{X}={\frac {1} {n}}\rho^{-1}(|-nK_{Y}|)$ for
some  natural number $n\gg 0$. Then $K_{X}+M_{X}\sim_\mathbb{Q} 0$
and the set of centers of canonical singularities $\mathbb{CS}(X,
M_{X})$ is not empty.
\end{theorem}

\begin{proof}
See \cite{Ch00a}, \cite{Ch03b} and \cite{ChPa04}.
\end{proof}

\begin{theorem}
\label{theorem:Shokurov-vanishing} Let $(X, B_{X})$ be a log pair
with effective $B_{X}$, $\mathcal{I}(X, B_{X})$ be an ideal sheaf
of the log canonical singularities subscheme $\mathcal{L}(X,
B_{X})$, and let $H$ be a nef and big divisor on $X$ such that
$K_{X}+B_{X}+H$ is a Cartier divisor. Then $H^{i}(X,
\mathcal{I}(X, B_{X})\otimes (K_{X}+B_{X}+H))=0$ for $i>0$.
\end{theorem}

\begin{proof}
See \cite{Sh92}, \cite{Ko91}, \cite{Ko97}, \cite{Am99} or
\cite{Ch03b}.
\end{proof}

\begin{theorem}
\label{theorem:log-adjunction} Let $(X, B_{X})$ be a log pair,
$B_{X}$ be a effective boundary such that $\lfloor
B_{X}\rfloor=\emptyset$, and let $S\subset X$ be an effective
irreducible divisor such that the divisor $K_{X}+S+B_{X}$ is
$\mathbb{Q}$-Cartier. Then $(X, S+B_{X})$ is purely log terminal
if and only if $(S, \mathrm{Diff}_{S}(B_{X}))$ is Kawamata log
terminal.
\end{theorem}

\begin{proof}
See Theorem~17.6 in \cite{Ko91} or Theorem~7.5 in \cite{Ko97}.

\end{proof}

\begin{corollary}
\label{corollary:log-adjunction} Let $(X, B_{X})$ be a log pair
with effective $B_{X}$, $H$ be an effective Cartier divisor on
$X$, $Z\in \mathbb{CS}(X, B_{X})$, both $X$ and $H$ are smooth in
the generic point of $Z\subset H\not\subset\mathrm{Supp}(B_{X})$.
Then the set of centers of log canonical singularities
$\mathbb{LCS}(H, B_{X}\vert_{H})$ is not empty.
\end{corollary}

\begin{theorem}
\label{theorem:Iskovskikh} Let $X$ be a smooth variety,
$\mathrm{dim}(X)\ge 3$, $M_{X}$ be an effective movable boundary
on the variety $X$, and the set $\mathbb{CS}(X, M_{X})$ contains a
closed point $O\in X$. Then the inequality
$\mathrm{mult}_{O}(M_{X}^{2})\geq 4$ holds and the equality
implies $\mathrm{mult}_{O}(M_{X})=2$ and $\mathrm{dim}(X)=3$.
\end{theorem}

\begin{proof}
See \cite{IsMa71}, \cite{Co95}, \cite{Pu00}, \cite{Co00} and
\cite{Kaw01}.
\end{proof}

\begin{theorem}
\label{theorem:Corti} Let $X$ be a variety, $\mathrm{dim}(X)\ge
3$, and $B_{X}$ be an effective boundary on $X$ such that the set
$\mathbb{CS}(X, B_{X})$ contains an ordinary double point $O$ of
$X$. Then $\mathrm{mult}_{O}(B_{X})\geq 1$ and the equality
implies $\mathrm{dim}(X)=3$.
\end{theorem}

\begin{proof}
The claim is implied by Theorem~3.10 in \cite{Co00} and
Theorem~\ref{theorem:log-adjunction}.
\end{proof}

\begin{proposition}
\label{proposition:double-covers} Let $\tau:V\to\mathbb{P}^{k}$ be
a double cover ramified in a smooth hypersurface
$S\subset\mathbb{P}^{k}$ of degree $2d$ such that $2\le d\le k-1$,
$B_{V}$ be an effective boundary on $V$ such that ${\frac {1}
{r}}B_{V}\sim_{\mathbb{Q}}
\tau^{*}(\mathcal{O}_{\mathbb{P}^{k}}(1))$ for some positive
rational number $r<1$. Then the set of centers of log canonical
singularities $\mathbb{LCS}(V, B_{V})$ is empty.
\end{proposition}

\begin{proof}
Let $C\subset V$ be an irreducible curve such that $\tau(C)\subset
S$ and the inequality $\mathrm{mult}_{C}(B_{V})\ge 1$ holds. Take
a point $O$ on the curve $\tau(C)$ and a hyperplane
$\Pi\subset\mathbb{P}^{k}$ that tangents $S$ at the point $O$. Fix
a line $L\subset \Pi$ passing through $O$. Let ${\hat
L}=\tau^{-1}(L)$. Then ${\hat L}$ is singular at ${\hat
O}=\tau^{-1}(O)$ and a component of ${\hat L}$ is contained in
$\mathrm{Supp}(B_{V})$, because otherwise
$$
2>2r={\hat L}\cdot B_{V}\ge\mathrm{mult}_{\hat O}(\hat L)\mathrm{mult}_{C}(B_{V})\ge 2%
$$
which is a contradiction. On the other hand, $\Pi$ tangents $S$ in
finitely many points (see \cite{FuLa81}, \cite{Ish82},
\cite{Pu95}, \cite{Zak93}). Hence, the curve ${\hat L}$ spans $V$
when we vary the point $O$ on the curve $\tau(C)$ and the line
$L\subset\Pi$. The latter is a contradiction, because ${\hat
L}\subset\mathrm{Supp}(B_{V})$.

Suppose that $\mathbb{LCS}(V, B_{V})$ contains a subvariety
$Z\subset V$ of dimension at least two. Then
$\mathrm{mult}_{Z}(B_{V})\ge 1$ and the set $Z\cap \tau^{-1}(S)$
contains some curve ${\hat C}\subset V$. Then
$\mathrm{mult}_{{\hat C}}(B_{V})\ge 1$ and $\tau({\hat C})\subset
S$, but we already prove that this is impossible. Hence, the set
$\mathbb{LCS}(V, B_{V})$ does not contains subvarieties of
dimension at least two.

Suppose that the set $\mathbb{LCS}(V, B_{V})$ contains a curve on
$V$. Consider a union $T\subset V$ of all curves in the set
$\mathbb{LCS}(V, B_{V})$. We may consider $T$ as a possibly
reducible curve on $V$. Let $Y$ be a sufficiently general divisor
in the linear system
$|\tau^{*}(\mathcal{O}_{\mathbb{P}^{k}}(1))|$,
$\gamma=\tau\vert_{Y}$ and $B_{Y}=B_{V}\vert_{Y}$. Then the
variety $Y$ is smooth, $Y\not\subset\mathrm{Supp}(B_{V})$, and
$\gamma:Y\to \mathbb{P}^{k-1}$ is a double cover branched over a
smooth hypersurface of degree $2d$. The generality in the choice
of $Y$ implies that the set $\mathbb{LCS}(Y, B_{Y})$ does not
contain subvarieties of $Y$ of positive dimension. Moreover, the
set $\mathbb{LCS}(Y, B_{Y})$ is not empty, i.e. it contains all
points of $T\cap Y$. Consider a Cartier divisor
$$
F=K_{Y}+B_{Y}+(1-r)H\sim (d-k-1)H
$$
where $H=\gamma^{*}(\mathcal{O}_{\mathbb{P}^{k-1}}(1))$. The
sequence of groups
$$
H^{0}(\mathcal{O}_{Y}(F))\to H^{0}(\mathcal{O}_{\mathcal{L}(Y,
B_{Y})}(F))\to 0
$$
is exact by Theorem~\ref{theorem:Shokurov-vanishing} where
$\mathcal{L}(Y, B_{Y})$ is a log canonical singularities subscheme
of $(Y, B_{Y})$. On the other hand, $\mathrm{Supp}(\mathcal{L}(Y,
B_{Y}))$ consists of all points in $T\cap Y$. Hence,
$H^{0}(\mathcal{O}_{\mathcal{L}(Y,
B_{Y})}(F))=H^{0}(\mathcal{O}_{\mathcal{L}(Y, B_{Y})})$. The
latter contradicts $d<k-1$, because $H^{0}(\mathcal{O}_{Y}(F))=0$
for $d<k-1$. However, in the case when $d=k-1$ the latter implies
that the set $T\cap Y$ consists of a single point, because
$H^{0}(\mathcal{O}_{Y}(F))=\mathbb{C}$ for $d=k-1$.

We proved that the assumption that the set $\mathbb{LCS}(V,
B_{V})$ contains some curve on $V$ implies that $d=k-1$, the set
$\mathbb{LCS}(V, B_{V})$ contains a single curve ${\bar C}\subset
V$ such that $\tau(\bar C)\subset\mathbb{P}^{k}$ is a line,
$\tau\vert_{\bar C}$ is an isomorphism, and the inequality
$\mathrm{mult}_{\bar C}(B_{V})\ge 1$ holds. On the other hand, we
already proved that the latter implies $\tau({\bar C})\not\subset
S$. Therefore, there is an irreducible reduced curve $\tilde
C\subset V$ such that ${\bar C}\ne {\tilde C}$ and $\tau({\bar
C})=\tau({\tilde C})$.

Let $D_{1},\ldots,D_{k-2}$ be sufficiently general divisors in
$|\tau^{*}(\mathcal{O}_{\mathbb{P}^{k}}(1))|$ passing through the
curves ${\bar C}$ and ${\tilde C}$. Put $D=\cap_{i=1}^{k-2}D_{i}$.
Then $D$ is a smooth surface, and both ${\bar C}$ and ${\tilde C}$
are smooth rational curves on $D$. By the adjunction formula the
self-intersections of the curves ${\bar C}$ and ${\tilde C}$ on
the surface $D$ are equal to $1-d$. Therefore, ${\bar
C}^{2}={\tilde C}^{2}<0$ due to $d>2$.

By construction we have $D\not\subset\mathrm{Supp}(B_{V})$.
Therefore, we can consider a boundary $B_{D}=B_{V}\vert_{D}$. The
generality in the choice of $D$ implies
$$
B_{D}=\mathrm{mult}_{{\bar C}}(B_{V}){\bar C}+\mathrm{mult}_{{\tilde C}}(B_{V}){\tilde C}+\Delta%
$$
where $\Delta$ is an effective divisor on the surface $D$ such
that $\mathrm{Supp}(\Delta)$ does not contain both curves ${\bar
C}$ and ${\tilde C}$. On the other hand, the equivalence
$$
B_{D}\sim_{\mathbb{Q}} r({\bar C}+{\tilde C})
$$
holds. In particular, the equivalence
$$
(r-\mathrm{mult}_{{\tilde C}}(B_{V})){\tilde C}\sim_{\mathbb{Q}}(\mathrm{mult}_{{\bar C}}(B_{V})-r){\bar C}+\Delta%
$$
holds. Therefore, $\mathrm{mult}_{{\tilde C}}(B_{V})\ge r$ due to
${\tilde C}^{2}<0$. Thus, the equivalence
$$
-\Delta\sim_{\mathbb{Q}}(\mathrm{mult}_{{\bar C}}(B_{V})-r){\bar C}+(\mathrm{mult}_{{\tilde C}}(B_{V})-r){\tilde C}%
$$
implies $\Delta=\emptyset$ and $\mathrm{mult}_{{\tilde
C}}(B_{V})=\mathrm{mult}_{{\bar C}}(B_{V})=r$. The latter is
impossible, because $\mathrm{mult}_{{\bar C}}(B_{V})\ge 1$ and
$r<1$. Therefore, the set $\mathbb{LCS}(V, B_{V})$ does not
contain subvarieties of positive dimension.

Suppose that $\mathbb{LCS}(V, B_{V})$ contains a closed point $O$
on $V$. Let
$$
E=K_{V}+B_{V}+(1-r)H
$$
where $H=\tau^{*}(\mathcal{O}_{\mathbb{P}^{k}}(1))$. Then $E$ is a
Cartier divisor and $H^{0}(\mathcal{O}_{V}(E))=0$, because $E\sim
(d-k)H$ and $d\le k-1$. However, the sequence
$$
H^{0}(\mathcal{O}_{V}(E))\to H^{0}(\mathcal{O}_{\mathcal{L}(V,
B_{V})}(E))\to 0
$$
is exact by Theorem~\ref{theorem:Shokurov-vanishing} where
$\mathcal{L}(V, B_{V})$ is a log canonical singularities subscheme
of $(V, B_{V})$. On the other hand, $\mathrm{Supp}(\mathcal{L}(V,
B_{V})$ consists of finite number of points of $V$. Hence,
$H^{0}(\mathcal{O}_{\mathcal{L}(V,
B_{V})}(E))=H^{0}(\mathcal{O}_{\mathcal{L}(V, B_{V})})$ which is a
contradiction. Thus, the set $\mathbb{LCS}(V, B_{V})$ is empty.

\end{proof}

\begin{proposition}
\label{proposition:projective-space} Let $S\subset\mathbb{P}^{n}$
be a smooth hypersurface of degree $d\ge 2$ and $B$ be an
effective boundary on $\mathbb{P}^{n}$ such that ${\frac {1}
{r}}B\sim_{\mathbb{Q}} \mathcal{O}_{\mathbb{P}^{n}}(1)$ for a
positive rational number $r<1$. Then the set of centers of log
canonical singularities $\mathbb{LCS}(\mathbb{P}^{n},
B+{\frac{1}{2}}S)$ is empty for $d\le 2(n-1)$.
\end{proposition}

\begin{proof}
Let $Z\in \mathbb{LCS}(\mathbb{P}^{n}, B+{\frac{1}{2}}S)$ be a
center of maximal dimension. Then
$$
r+{\frac{1}{2}}\ge \mathrm{mult}_{Z}(B)+{\frac{1}{2}}\mathrm{mult}_{Z}(S)\ge \mathrm{mult}_{Z}(B+{\frac{1}{2}}S)\ge 1%
$$
which implies $Z\subset S$ and $Z\ne S$. Hence,
$\mathrm{dim}(Z)<n-1$.

Suppose that $Z$ is a closed point. Let
$$
E=K_{\mathbb{P}^{n}}+B+{\frac{1}{2}}S+(n-{\frac{d}{2}}-r)H
$$
where $H\sim\mathcal{O}_{\mathbb{P}^{n}}(1)$. Then $E$ is a
Cartier divisor, $n-{\frac{d}{2}}-r>0$ and the equivalence $E\sim
-H$ holds. Thus, $H^{0}(\mathcal{O}_{\mathbb{P}^{n}}(E))=0$. The
sequence
$$
H^{0}(\mathcal{O}_{\mathbb{P}^{n}}(E))\to
H^{0}(\mathcal{O}_{\mathcal{L}(\mathbb{P}^{n},
B+{\frac{1}{2}}S)}(E))\to 0
$$
is exact by Theorem~\ref{theorem:Shokurov-vanishing} where
$\mathcal{L}(\mathbb{P}^{n}, B+{\frac{1}{2}}S)$ is a log canonical
singularities subscheme of $(\mathbb{P}^{n}, B+{\frac{1}{2}}S)$.
However, $\mathrm{Supp}(\mathcal{L}(\mathbb{P}^{n},
B+{\frac{1}{2}}S))$ consists of finite number of closed points of
$\mathbb{P}^{n}$. Hence,
$$
H^{0}(\mathcal{O}_{\mathcal{L}(\mathbb{P}^{n}, B+{\frac{1}{2}}S)}(E))=H^{0}(\mathcal{O}_{\mathcal{L}(\mathbb{P}^{n}, B+{\frac{1}{2}}S)})%
$$
which is a contradiction. Thus, $\mathrm{dim}(Z)>0$.

Rewrite $B+{\frac{1}{2}}S$ as $D+\lambda S$ for an effective
boundary $D$ on $\mathbb{P}^{n}$ and a positive rational $\lambda$
such that $S\not\subset\mathrm{Supp}(D)$. Then $\lambda<1$ and
$D\sim_{\mathbb{Q}} \mu H$ for a positive rational number $\mu<1$.
In particular, $Z\subset S$ is a center of log canonical
singularities of log pair $(\mathbb{P}^{n}, D+S)$. Thus,
Theorem~\ref{theorem:log-adjunction} implies $\mathbb{LCS}(S,
D\vert_{S})\ne\emptyset$. Moreover,
Theorem~\ref{theorem:log-adjunction} implies the existence of a
subvariety $T\subset S$ such that $T\in \mathbb{LCS}(S,
D\vert_{S})$ and $Z\subseteq T$. In particular, the inequalities
$\mathrm{dim}(T)\ge 1$ and $\mathrm{mult}_{T}(D\vert_{S})\ge 1$
hold, where $S$ is smooth by assumption. The latter is impossible
due to \cite{Pu95}. Namely, let $C$ be a curve in $T$,
$Y\subset\mathbb{P}^{n}$ be a general cone over $C$ and ${\tilde
C}\subset S$ be a residual curve to the curve $C$ defined as
$C\cup {\tilde C}=Y\cap S$. Then $\mathrm{mult}_{C}(D\vert_{S})\ge
1$, the intersection $C\cap {\tilde C}$ consists of
$(\mathrm{deg}(S)-1)\mathrm{deg}(C)$ different points in a
set-theoretic sense, and ${\tilde C}\not\subset\mathrm{Supp}(D)$.
In particular,
$$
\mathrm{deg}(D\vert_{{\tilde C}})\ge (\mathrm{deg}(S)-1)\mathrm{deg}(C)\mathrm{mult}_{C}(D\vert_{S})\ge (\mathrm{deg}(S)-1)\mathrm{deg}(C),%
$$
but $\mathrm{deg}(D\vert_{{\tilde
C}})=\mu(\mathrm{deg}(S)-1)\mathrm{deg}(C)$, which is a
contradiction.
\end{proof}

\section{The proof of Theorem~\ref{theorem:main}.}
\label{section:proof-of-main-result}

Let $\pi:X\to \mathbb{P}^{n}$ be a double cover ramified in a
hypersurface $F\subset\mathbb{P}^{n}$ of degree $2n$ with isolated
singularities such that $n\ge 4$ and every singular point $O$ of
$F$ is an ordinary singular point and $\mathrm{mult}_{O}(F)\le
2(n-2)$.

\begin{lemma}
\label{lemma:factoriality} The variety $V$ is a Fano variety with
terminal $\mathbb{Q}$-factorial singularities such that
$\mathrm{Cl}(X)\cong\mathrm{Pic}(X)\cong\mathbb{Z}$.
\end{lemma}

\begin{proof}
The ampleness of the divisor $-K_{X}$ and the terminality of $X$
are obvious. Consider a Weil divisor $D$ on the variety $X$. To
prove the claim it is enough to show that $D\sim
\pi^{*}(\mathcal{O}_{\mathbb{P}^{n}}(r))$ for some
$r\in\mathbb{Z}$.

Let $H$ be a general divisor in
$|\pi^{*}(\mathcal{O}_{\mathbb{P}^{n}}(k))|$ for $k\gg 0$. Then
$H$ is a smooth complete intersection in $\mathbb{P}(1^{n+1}, n)$
and $\mathrm{dim}(X)\ge 3$. Therefore, the group $\mathrm{Pic}(H)$
is generated by
$\pi^{*}(\mathcal{O}_{\mathbb{P}^{n}}(1))\vert_{H}$ by
Th\'eor\'eme 3.13 of Exp. XI in \cite{Gro65} (see Lemma 3.2.2 in
\cite{Do82}, Lemma 3.5 in \cite{CPR} or \cite{CalLy94}). Thus,
there is an integer $r$ such that $D\vert_{H}\sim
\pi^{*}(\mathcal{O}_{\mathbb{P}^{n}}(r))\vert_{H}$.

Let $\mathrm{\Delta}=D-\pi^{*}(\mathcal{O}_{\mathbb{P}^{n}}(r))$.
The sequence of sheaves
$$
0\to\mathcal{O}_{X}(\mathrm{\Delta})\otimes \pi^{*}(\mathcal{O}_{\mathbb{P}^{n}}(-k))\to\mathcal{O}_{X}(\mathrm{\Delta})\to\mathcal{O}_{H}\to 0%
$$
is exact, because the sheaf $\mathcal{O}_{X}(\mathrm{\Delta})$ is
locally free in the neighborhood of the divisor $H$. Therefore,
the sequence of groups
$$
0\to H^{0}(\mathcal{O}_{X}(\mathrm{\Delta}))\to H^{0}(\mathcal{O}_{H})\to H^{1}(\mathcal{O}_{X}(\mathrm{\Delta})\otimes \pi^{*}(\mathcal{O}_{\mathbb{P}^{n}}(-k)))%
$$
is exact. On the other hand, there is an exact sequence of sheaves
$$
0\to\mathcal{O}_{X}(\mathrm{\Delta})\to\mathcal{E}\to\mathcal{F}\to 0%
$$
where $\mathcal{E}$ is a locally free sheaf and $\mathcal{F}$ is a
torsion free sheaf, because the sheaf
$\mathcal{O}_{X}(\mathrm{\Delta})$ is reflexive (see
\cite{Har80}). Hence, the sequence of groups
$$
H^{0}(\mathcal{F}\otimes\mathcal{O}_{X}(-H))\to H^{1}(\mathcal{O}_{X}(\mathrm{\Delta}-H))\to H^{1}(\mathcal{E}\otimes\mathcal{O}_{X}(-H))%
$$
is exact. However,
$H^{0}(\mathcal{F}\otimes\mathcal{O}_{X}(-H))=0$, because the
sheaf $\mathcal{F}$ has no torsion, and
$H^{1}(\mathcal{E}\otimes\mathcal{O}_{X}(-H))=0$ by the lemma of
Enriques-Severi-Za\-riski (see \cite{Za52}). Thus, we have
$$
H^{1}(\mathcal{O}_{X}(\mathrm{\Delta})\otimes\pi^{*}(\mathcal{O}_{\mathbb{P}^{n}}(-k)))=0%
$$
and $H^{0}(\mathcal{O}_{X}(\mathrm{\Delta}))=\mathbb{C}$. The same
method gives
$H^{0}(\mathcal{O}_{X}(-\mathrm{\Delta}))=\mathbb{C}$, i.e. the
divisor $\mathrm{\Delta}$ is rationally equivalent to zero.
\end{proof}

Suppose that $X$ is not birationally super-rigid. Then there is a
movable log pair $(X, M_{X})$ such that $M_{X}$ is effective, the
set of centers of canonical singularities $\mathbb{CS}(X,M_{X})$
is not empty and the divisor $-(K_{X}+M_{X})$ is ample by
Theorem~\ref{theorem:Nother-Fano}. Let $Z$ be an element of the
set $\mathbb{CS}(X,M_{X})$.

\begin{lemma}
\label{lemma:smooth-points} The subvariety $Z\subset X$ is not a
smooth point of $X$.
\end{lemma}

\begin{proof}
Let $Z$ be a smooth point of $X$. Then
$\mathrm{mult}_{Z}(M_{X}^{2})>4$ by
Theorem~\ref{theorem:Iskovskikh}. Consider $n-2$ general divisors
$H_{1},\ldots,H_{n-2}$ in
$|\pi^{*}(\mathcal{O}_{\mathbb{P}^{n}}(1))|$ that pass through the
point $Z$. Then
$$
2>M^{2}_{X}\cdot H_{1}\cdots H_{n-2}\ge \mathrm{mult}_{Z}(M^{2}_{X})\mathrm{mult}_{Z}(H_{1})\cdots\mathrm{mult}_{Z}(H_{n-2})>4%
$$
which is a contradiction.
\end{proof}

\begin{lemma}
\label{lemma:singular-points} The subvariety $Z\subset X$ is not a
singular point of $X$.
\end{lemma}

\begin{proof}
The variety $X$ can be given as a hypersurface
$$
y^{2}=f_{2n}(x_{0},\ldots,x_{n})\subset\mathbb{P}(1^{n+1},n)\cong\mathrm{Proj}(\mathbb{C}[x_{0},\ldots,x_{n},y])%
$$
where $f_{2n}$ is a homogeneous polynomial of degree $2n$. Suppose
that $Z$ is a singular point of $X$. Then $O=\pi(Z)$ is an
ordinary singular point on the hypersurface
$F\subset\mathbb{P}^{n}$. There are two possible cases, i.e.
$\mathrm{mult}_{O}(F)$ is even or odd. We handle them separately.

Suppose $\mathrm{mult}_{O}(F)=2m\ge 2$ for some $m\in\mathbb{N}$.
By the initial assumption $m\le n-2$. There is a weighted blow up
$\beta:U\to\mathbb{P}(1^{n+1},n)$ of the point $Z$ with weights
$(m,1^{n})$ such that the proper transform $V\subset U$ of the
variety $X$ is non-singular in the neighborhood of the
$\beta$-exceptional divisor $E$. The morphism $\beta$ induces a
birational morphism $\alpha:V\to X$ with an exceptional divisor
$G\subset V$. Then $E\vert_{V}=G$ and $G$ is a smooth hypersurface
in $E\cong\mathbb{P}(1^{n},m)$ which can be given by
$$
z^{2}=g_{2m}(t_{1},\ldots,t_{n})\subset\mathbb{P}(1^{n},m)\cong\mathrm{Proj}(\mathbb{C}[t_{1},\ldots,t_{n},z])%
$$
where $g_{2m}$ is a homogeneous polynomial of degree $2m$.

Let $M_{V}=\alpha^{-1}(M_{X})$ and $\mathrm{mult}_{Z}(M_{X})$ be a
positive rational number such that $M_{V}\sim_{\mathbb{Q}}
\alpha^{*}(M_{X})-\mathrm{mult}_{Z}(M_{X})G$. Then the equivalence
$$
K_{V}+M_{V}\sim_{\mathbb{Q}}\alpha^{*}(K_{X}+M_{X})+(n-1-m-\mathrm{mult}_{Z}(M_{X}))G%
$$
holds. However, the linear system $|\alpha^{*}(-K_{X})-G|$ is free
and gives a fibration $\psi:V\to\mathbb{P}^{n-1}$ such that
$\psi=\chi\circ\pi\circ\alpha$ where
$\chi:\mathbb{P}^{n}\dashrightarrow \mathbb{P}^{n-1}$ is a
projection from $O$. Let $C$ be a general fiber of $\psi$. Then
$$
M_{V}\cdot C=2(1-\mathrm{mult}_{Z}(M_{X}))+\alpha^{*}(K_{X}+M_{X})\cdot C<2(1-\mathrm{mult}_{Z}(M_{X}))%
$$
because $-(K_{X}+M_{X})$ is ample. Thus,
$\mathrm{mult}_{Z}(M_{X})<1$. The latter contradicts
Theorem~\ref{theorem:Corti} in the case of $m=1$. Thus, $m>1$. On
the other hand, the inequality
$(n-1-m-\mathrm{mult}_{Z}(M_{X}))>0$ implies the existence of a
center $\Delta\in\mathbb{CS}(V, M_{V})$ such that $\Delta\subset
G$. Hence, $\mathbb{LCS}(G, M_{V}\vert_{G})\ne\emptyset$ by
Corollary~\ref{corollary:log-adjunction}. The latter contradicts
Proposition~\ref{proposition:double-covers}.

Therefore, $\mathrm{mult}_{O}(F)=2k+1\ge 3$ for $k\in\mathbb{N}$.
Then $k\le n-3$ by the initial assumption. Let
$\lambda:W\to\mathbb{P}^{n}$ be a blow up of $O$, $\Lambda$ be an
exceptional divisor of the birational morphism $\lambda$, and
${\tilde F}\subset W$ be a proper transform of the hypersurface
$F$. Then ${\tilde F}$ is smooth in the neighborhood of the
exceptional divisor $\Lambda$ and $S={\tilde F}\cap \Lambda\subset
\Lambda\cong\mathbb{P}^{n-1}$ is a smooth hypersurface of degree
$2k+1$. Let ${\tilde\pi}:{\tilde X}\to W$ be a double cover
ramified in the effective divisor
$$
{\tilde F}\cup\Lambda\sim 2(\lambda^{*}(\mathcal{O}_{\mathbb{P}^{n}}(n))-k\Lambda)%
$$
which is singular only in $S$. Then $W$ is smooth outside of
${\tilde S}={\tilde \pi}^{-1}(S)$ and the singularities of $W$
along $\tilde S$ is of type $\mathbb{A}_{1}\times
\mathbb{C}^{n-2}$, i.e. a two-dimensional ordinary double point
along $\tilde S$. Let $\Xi={\tilde \pi}^{-1}(\Lambda)$. Then
$\Xi\cong\mathbb{P}^{n-1}$ and there is a birational morphism
$\xi:{\tilde X}\to X$ contracting $\Xi$ to the point $Z$ such that
$\pi\circ\xi=\lambda\circ {\tilde \pi}$. The birational morphism
$\xi$ is a restriction of the weighted blow up of
$\mathbb{P}(1^{n+1},n)$ at $Z$ with weights $(2k+1,2^{n})$.

Let $M_{\tilde X}=\xi^{-1}(M_{X})$ and $\mathrm{mult}_{Z}(M_{X})$
be a positive rational number such that $M_{\tilde
X}\sim_{\mathbb{Q}} \xi^{*}(M_{X})-\mathrm{mult}_{Z}(M_{X})\Xi$.
Then the equivalence
$$
K_{\tilde X}+M_{\tilde X}\sim_{\mathbb{Q}}\xi^{*}(K_{X}+M_{X})+(2(n-1-k)-\mathrm{mult}_{Z}(M_{X}))\Xi%
$$
holds. On the other hand, the linear system
$|\xi^{*}(-K_{X})-2\Xi|$ is free and gives a fibration
$\omega:{\tilde X}\to\mathbb{P}^{n-1}$ such that
$\omega=\chi\circ\pi\circ\xi$, where $\chi$ is a projection of
$\mathbb{P}^{n}$ to $\mathbb{P}^{n-1}$ from $O$. Intersecting
$M_{\tilde X}$ with a general fiber of $\omega$ we get
$\mathrm{mult}_{Z}(M_{X})<2$. Thus,
$(2(n-1-k)-\mathrm{mult}_{Z}(M_{X}))>0$ which implies the
existence of a center
$$
Z\in\mathbb{CS}({\tilde X}, M_{\tilde X}-(2(n-1-k)-\mathrm{mult}_{Z}(M_{X}))\Xi)%
$$
such that $Z\subset G$. Hence,
$$
Z\in\mathbb{LCS}({\tilde X}, M_{\tilde X}-(2(n-1-k)-\mathrm{mult}_{Z}(M_{X}))\Xi+2\Xi)%
$$
because $2\Xi$ is a Cartier divisor. However,
$$
\mathbb{LCS}({\tilde X}, M_{\tilde X}-(2(n-2-k)-\mathrm{mult}_{Z}(M_{X}))\Xi)\subset\mathbb{LCS}({\tilde X}, M_{\tilde X}+\Xi)%
$$
due to $2k+1\le 2(n-2)$, which implies
$$
\mathbb{LCS}(\Xi, \mathrm{Diff}_{\Xi}(M_{\tilde X}))=\mathbb{LCS}(\Xi, M_{\tilde X}\vert_{\Xi}+\mathrm{Diff}_{\Xi}(0))\ne\emptyset%
$$
by Theorem~\ref{theorem:log-adjunction}. However,
$\mathrm{Diff}_{\Xi}(0)={\frac{1}{2}}{\tilde S}$ (see \cite{Ko91},
\cite{Pr01}) and
$$
M_{\tilde X}\vert_{\Xi}\sim_{\mathbb{Q}}-\mathrm{mult}_{Z}(M_{X})\Xi\vert_{\Xi}\sim_{\mathbb{Q}}{\frac{\mathrm{mult}_{Z}(M_{X})}{2}}H,%
$$
where $H$ is a hyperplane on $\Xi\cong\mathbb{P}^{n-1}$.
Therefore, the set of log canonical singularities
$\mathbb{LCS}(\Xi, M_{\tilde X}\vert_{\Xi}+{\frac{1}{2}}{\tilde
S})$ is empty by Proposition~\ref{proposition:projective-space},
which is a contradiction.
\end{proof}

\begin{lemma}
\label{lemma:codimension-big} The inequality
$\mathrm{codim}(Z\subset X)>2$ is impossible.
\end{lemma}

\begin{proof}
Suppose that $\mathrm{codim}(Z\subset X)>2$. Then
$\mathrm{dim}(Z)\ne 0$ by Lemmas~\ref{lemma:smooth-points} and
\ref{lemma:singular-points}. Thus,
$\mathrm{mult}_{Z}(M_{X}^{2})\ge 4$ by
Theorem~\ref{theorem:Iskovskikh}. Take a point $O$ on $Z$ and
sufficiently general divisors $H_{1},\ldots,H_{n-2}$ in
$|\pi^{*}(\mathcal{O}_{\mathbb{P}^{n}}(1))|$ that pass through the
point $O$. Then
$$
2>M^{2}_{X}\cdot H_{1}\cdots H_{n-2}\ge\mathrm{mult}_{Z}(M^{2}_{X})\ge 4%
$$
which is a contradiction.
\end{proof}

\begin{lemma}
\label{lemma:codimension-two} The equality
$\mathrm{codim}(Z\subset X)=2$ is impossible.
\end{lemma}

\begin{proof}
Suppose $\mathrm{codim}(Z\subset X)=2$. Then
$\mathrm{mult}_{Z}(M_{X})\ge 1$. Take sufficiently general
divisors $H_{1},\ldots,H_{n-2}$ in
$|\pi^{*}(\mathcal{O}_{\mathbb{P}^{n}}(1))|$. Then
$$
2>M^{2}_{X}\cdot H_{1}\cdots
H_{n-2}\ge\mathrm{mult}^{2}_{Z}(M_{X})Z\cdot H_{1}\cdots
H_{n-2}\ge Z\cdot H_{1}\cdots H_{n-2},%
$$
because $-(K_{X}+M_{X})$ is ample and $K_{X}\sim
\pi^{*}(\mathcal{O}_{\mathbb{P}^{n}}(-1))$. Thus, $\pi(Z)$ is a
linear subspace in $\mathbb{P}^{n}$ of dimension $n-2$ and
$\pi\vert_{Z}$ is an isomorphism.

Let $V=\cap_{i=1}^{n-3}H_{i}$, $C=Z\cap V$, $M_{V}=M_{X}\vert_{V}$
and $\tau=\pi\vert_{V}$. Then $V$ is a smooth 3-fold, $C\subset V$
is a curve, $M_{V}$ is movable, $\tau:V\to\mathbb{P}^{3}$ is a
double cover branched over a smooth hypersurface $S\subset
\mathbb{P}^{3}$ of degree $2n$,  $\tau(C)$ is a line in
$\mathbb{P}^{3}$, $\tau\vert_{C}$ is an isomorphism. Moreover,
$\tau^{*}(\mathcal{O}_{\mathbb{P}^{3}}(1))-M_{V}$ is an ample
divisor and $\mathrm{mult}_{C}(M_{V})=\mathrm{mult}_{Z}(M_{X})$.

Suppose $\tau(C)\not\subset S$. Then there is an irreducible curve
$\tilde C\subset V$ such that $C\ne {\tilde C}$ and
$\tau(C)=\tau({\tilde C})$. Take a general divisor $D\in
|\tau^{*}(\mathcal{O}_{\mathbb{P}^{3}}(1))|$ passing through $C$.
Then $D$ is a smooth surface, $C$ and ${\tilde C}$ are smooth
rational curves. By the adjunction formula $C^{2}={\tilde
C}^{2}=1-n<0$ on the surface $D$. Consider a boundary
$M_{D}=M_{V}\vert_{D}$. The boundary $M_{D}$ is no longer movable.
However, the generality in the choice of $D$ implies
$$
M_{D}=\mathrm{mult}_{C}(M_{V})C+\mathrm{mult}_{{\tilde C}}(M_{V}){\tilde C}+\Delta%
$$
where $\Delta$ is a movable boundary on $D$. However,
$M_{V}\sim_{\mathbb{Q}} rD$ for some rational number $r<1$. Hence,
the equivalence
$$
(r-\mathrm{mult}_{{\tilde C}}(M_{V})){\tilde C}\sim_{\mathbb{Q}}(\mathrm{mult}_{C}(M_{V})-r)C+\Delta%
$$
holds. The inequality ${\tilde C}^{2}<0$ implies
$\mathrm{mult}_{{\tilde C}}(M_{V})\ge r$. Let $H$ be a
sufficiently general divisor in the linear system
$|\tau^{*}(\mathcal{O}_{\mathbb{P}^{3}}(1))|$. Then
$$
2r^{2}=M^{2}_{V}\cdot H\ge
\mathrm{mult}^{2}_{C}(M_{V})+\mathrm{mult}^{2}_{{\tilde
C}}(M_{V})\ge 1+r^{2}
$$
which contradicts the inequality $r<1$.

Suppose that $\tau(C)\subset S$. Let $O$ be a general point on
$\tau(C)$ and $T$ be a hyperplane in $\mathbb{P}^{3}$ that
tangents $S$ at the point $O$. Consider a sufficiently general
line $L\subset T$ passing through $O$. Let ${\hat
L}=\tau^{-1}(L)$. Then ${\hat L}$ is singular at the point ${\hat
O}=\tau^{-1}(O)$. Therefore, ${\hat
L}\subset\mathrm{Supp}(M_{V})$, because otherwise
$$
2>{\hat L}\cdot M_{V}\ge\mathrm{mult}_{\hat O}(\hat L)\mathrm{mult}_{C}(M_{V})\ge 2%
$$
which is a contradiction. On the other hand, the curve ${\hat L}$
spans a divisor in the variety $V$ when we vary the line $L\subset
T$. The latter contradicts the movability of the boundary $M_{V}$.
\end{proof}

Therefore, Theorem~\ref{theorem:main} is proved.

\section{The proof of Theorems~\ref{theorem:second} and
\ref{theorem:third}.}
\label{section:proof-of-second-result}

Let $\pi:X\to \mathbb{P}^{n}$ be a double cover branched over an
hypersurface $F$ of degree $2n$ with isolated singularities,
$n=\mathrm{dim}(X)\ge 4$ and every singular point $O$ of the
hypersurface $F$ is an ordinary singular point of multiplicity
$\mathrm{mult}_{O}(F)\le 2(n-2)$. Let $\rho:X\dasharrow Y$ be a
birational map and $\tau:Y\to Z$ be an elliptic fibration. Take a
very ample divisor $H$ on the variety $Z$ and consider a linear
system $\mathcal{M}=\rho^{-1}(|\pi^{*}(H)|)$.

\begin{remark}
\label{remark:not-composed-from-a-pencil} The linear system
$\mathcal{M}$ is not composed from a pencil.
\end{remark}

Due to Lemma~\ref{lemma:factoriality} there is a positive rational
number $r$ such that the equivalence
$K_{X}+r\mathcal{M}\sim_\mathbb{Q} 0$ holds. Let
$M_{X}=r\mathcal{M}$. Then $\mathbb{CS}(X, M_{X})\ne\emptyset$ by
Theorem~\ref{theorem:elliptic-fibrations}. Let $Z$ be an element
of the set $\mathbb{CS}(X,M_{X})$.

\begin{lemma}
\label{lemma:smooth-points-II} The subvariety $Z\subset X$ is not
a smooth point of $X$.
\end{lemma}
\begin{proof}
See the proof of Lemma~\ref{lemma:smooth-points}.
\end{proof}

\begin{lemma}
\label{lemma:singular-points-II} Let $Z$ be a singular point of
$X$. Then $\mathrm{mult}_{O}(F)=2(n-2)$ and
$\tau\circ\rho=\gamma\circ\beta\circ\pi$, where $O=\pi(Z)$,
$\beta:\mathbb{P}^{n}\dasharrow\mathbb{P}^{n-1}$ is a projection
from the point $O$, and $\gamma:\mathbb{P}^{n-1}\dasharrow Y$ is a
birational map.
\end{lemma}

\begin{proof}
The point $O$ is an ordinary singular point of
$F\subset\mathbb{P}^{n}$ such that the inequality
$\mathrm{mult}_{O}(F)\le 2(n-2)$ holds. Suppose that the
multiplicity of the hypersurface $F$ at the point $O$ is even,
i.e. $\mathrm{mult}_{O}(F)=2m\ge 2$ for $m\in\mathbb{N}$. The
variety $X$ is a hypersurface in $\mathbb{P}(1^{n+1},n)$ of degree
$2n$, and there is a weighted blow up
$\beta:U\to\mathbb{P}(1^{n+1},n)$ of the point $Z$ with weights
$(m,1^{n})$ such that the proper transform $V\subset U$ of $X$ is
smooth near the exceptional divisor $E$ of $\beta$. The birational
morphism $\beta$ induces the birational morphism $\alpha:V\to X$.
Let $G$ be an exceptional divisor of $\alpha$. Then $E\vert_{V}=G$
and $G$ is a double cover of $\mathbb{P}^{n-1}$ branched over a
smooth hypersurface of degree $2m$.

Let $M_{V}=\alpha^{-1}(M_{X})$ and  $\mathrm{mult}_{Z}(M_{X})$ be
a positive rational number such that $M_{V}\sim_{\mathbb{Q}}
\alpha^{*}(M_{X})-\mathrm{mult}_{Z}(M_{X})G$. Then the equivalence
$$
K_{V}+M_{V}\sim_{\mathbb{Q}}\alpha^{*}(K_{X}+M_{X})+(n-1-m-\mathrm{mult}_{Z}(M_{X}))G%
$$
holds. On the other hand, the linear sister
$|\alpha^{*}(-K_{X})-G|$ is free and gives a fibration
$\psi:V\to\mathbb{P}^{n-1}$ such that
$\psi=\chi\circ\pi\circ\alpha$ where
$\chi:\mathbb{P}^{n}\dashrightarrow \mathbb{P}^{n-1}$ is a
projection from the point $O$. Let $C$ be a general fiber of
$\psi$. Then
$$
M_{V}\cdot C=2(1-\mathrm{mult}_{Z}(M_{Z}))
$$
and $g(C)=n-m+1$. Thus, $\mathrm{mult}_{Z}(M_{X})\le 1$. On the
other hand, the equality $\mathrm{mult}_{Z}(M_{X})=1$ implies that
$\psi$ and $\tau$ are birationally equivalent fibrations, i.e.
there is a birational map that maps the generic fiber of $\psi$
into the generic fiber of $\tau$. The latter is impossible in the
case of $m<n-2$, because $g(C)\ne 1$. In the case of $m=n-2$ the
equivalence of $\tau$ and $\psi$ implies the claim of the lemma.
Thus, we may assume $\mathrm{mult}_{Z}(M_{X})<1$ and proceed as in
the proof of Lemma~\ref{lemma:singular-points} to get a
contradiction.

Hence, we may assume that the multiplicity of the hypersurface $F$
at the point $O$ is odd. In this case the arguments above together
with the proof of Lemma~\ref{lemma:singular-points} give a
contradiction.
\end{proof}

\begin{lemma}
\label{lemma:codimension-big-II} The inequality
$\mathrm{codim}(Z\subset X)>2$ is impossible.
\end{lemma}

\begin{proof}
See the proof of Lemma~\ref{lemma:codimension-big}.
\end{proof}

\begin{lemma}
\label{lemma:codimension-two-II} The equality
$\mathrm{codim}(Z\subset X)=2$ is impossible.
\end{lemma}

\begin{proof}
Suppose $\mathrm{codim}(Z\subset X)=2$. Then
$\mathrm{mult}_{Z}(M_{X})\ge 1$. Take sufficiently general
divisors $H_{1},\ldots,H_{n-2}$ in
$|\pi^{*}(\mathcal{O}_{\mathbb{P}^{n}}(1))|$. Then
$$
2=M^{2}_{X}\cdot H_{1}\cdots
H_{n-2}\ge\mathrm{mult}^{2}_{Z}(M_{X})Z\cdot H_{1}\cdots
H_{n-2}\ge Z\cdot H_{1}\cdots H_{n-2},%
$$
and $k=Z\cdot H_{1}\cdots H_{n-2}$ is either $1$ or $2$.

Suppose $k=2$. Then for any two different divisors $D_{1}$ and
$D_{2}$ in the linear system $\mathcal{M}$ the intersection
$D_{1}\cap D_{2}$ coincide with $Z$ in the set-theoretic sense.
Let $p\not\in Z$ be a sufficiently general point and
$\mathcal{D}\subset \mathcal{M}$ be a linear subsystem of divisors
passing through the point $p$. Then $\mathcal{D}$ has no base
components, because $\mathcal{M}$ is not composed from a pencil.
Suppose that the divisors $D_{1}$ and $D_{2}$ are from
$\mathcal{D}$. Then in the set-theoretic sense
$$p\in D_{1}\cap D_{1}=Z$$
which is a contradiction. Therefore, $k=1$, i.e.
$\pi(Z)\subset\mathbb{P}^{n}$ is a linear subspace in of dimension
$n-2$ and $\pi\vert_{Z}$ is an isomorphism.

Suppose $\pi(Z)\not\subset F$. There is a subvariety ${\tilde
Z}\subset X$ of codimension two, such that $\pi({\tilde
Z})=\pi(Z)$ and ${\tilde Z}\ne Z$. The proof of
Lemma~\ref{lemma:codimension-two-II} gives
$$
\mathrm{mult}_{\tilde Z}(M_{X})=\mathrm{mult}_{Z}(M_{X})=1
$$
which leads to a contradiction as in the case of $k=2$. Thus,
$\pi(Z)\subset F$.

Consider a smooth 3-fold $V=\cap_{i=1}^{n-3}H_{i}$, a curve
$C=Z\cap V$, a movable boundary $M_{V}=M_{X}\vert_{V}$, a linear
system $\mathcal{D}=\mathcal{M}\vert_{V}$ that has no base
components, and a morphism $\tau=\pi\vert_{V}$. Then
$\tau:V\to\mathbb{P}^{3}$ is a double cover branched over a smooth
hypersurface $S\subset \mathbb{P}^{3}$ of degree $2n$,
$\tau(C)\subset S$ is a line, and $\tau\vert_{C}$ is an
isomorphism. Moreover, the equivalence
$$
M_{V}\sim_{\mathbb{Q}}\tau^{*}(\mathcal{O}_{\mathbb{P}^{3}}(1))
$$
holds and $\mathrm{mult}_{C}(M_{V})=\mathrm{mult}_{Z}(M_{X})\ge
1$.

Let $O$ be a general point on $\tau(C)$ and $T$ be a hyperplane in
$\mathbb{P}^{3}$ that tangents the hypersurface $S$ at the point
$O$. Consider a line $L\subset T$ passing through the point $O$.
Let ${\hat L}=\tau^{-1}(L)$. Then the curve ${\hat L}$ is singular
at the point ${\hat O}=\tau^{-1}(O)$. Therefore,
$\mathrm{mult}_{C}(M_{V})=1$, because
$$
2={\hat L}\cdot M_{V}\ge\mathrm{mult}_{\hat O}(\hat L)\mathrm{mult}_{C}(M_{V})\ge 2%
$$
and ${\hat L}$ spans a divisor when we vary the line $L\subset T$.
Let $f:U\to V$ be a blow up of $C$, $G$ be a $g$-exceptional
divisor, $M_{U}=f^{-1}(M_{V})$, $D$ be a general divisor in
$|(\tau\circ f)^{*}(\mathcal{O}_{\mathbb{P}^{3}}(1))-G|$,
$M_{D}=M_{U}\vert_{D}$. Then $D$ is smooth,
$$
M_{D}=\mathrm{mult}_{{\tilde C}}(M_{U}){\tilde C}+\Delta%
$$
where ${\tilde C}\subset G$ is a base curve of $|(\tau\circ
f)^{*}(\mathcal{O}_{\mathbb{P}^{3}}(1))-G|$ and $\Delta$ is a
movable boundary on $D$. The curve $\tilde C\subset G$ is a smooth
rational curve which dominates the curve $C$. By the adjunction
formula ${\tilde C}^{2}=1-n$ on the surface $D$. On the other
hand, the equivalence
$$
M_{D}\sim_{\mathbb{Q}}{\tilde C}%
$$
holds on $D$. The latter implies $\mathrm{mult}_{{\tilde
C}}(M_{U})=1$ and $\Delta=\emptyset$. After blowing up the curve
$\tilde C$ we see that the linear system $\mathcal{D}$ lies in the
fibers of the rational map given by the pencil $|(\tau\circ
f)^{*}(\mathcal{O}_{\mathbb{P}^{3}}(1))-G|$, which is impossible
because $\mathcal{D}$ is not composed from a pencil.
\end{proof}

Therefore, Theorem~\ref{theorem:second} is proved. The proof of
Theorem~\ref{theorem:third} is almost identical to the proof of
Theorem~\ref{theorem:second}. The only difference is that one must
use Theorem~\ref{theorem:canonical-Fanos} instead of
Theorem~\ref{theorem:elliptic-fibrations}.


\begin{thebibliography}{999}

\bibitem{Al91}
V.~Alexeev, \emph{General elephants of $\mathbb{Q}$-Fano 3-folds}, Comp. Math. \textbf{91} (1994), 91--116.%

\bibitem{Am99}
F.~Ambro, \emph{Ladders on Fano varieties}, J. Math. Sci. (New York) \textbf{94} (1999), 1126--1135.%

\bibitem{Ar62}
M.~Artin, \emph{Some numerical criteria of contractability of curves on algebraic surfaces}, Amer. J. Math. \textbf{84} (1962), 485--496.%

\bibitem{ArMu72}
M.~Artin, D.~Mumford, \emph{Some elementary examples of unirational varieties which are not rational}, Proc. London Math. Soc. \textbf{25} (1972), 75--95.%

\bibitem{Ba96}
W. Barth, \emph{Two projective surfaces with many nodes, admitting the symmetries of the icosahedron}, J. Algebraic Geometry \textbf{5} (1996), 173--186.%

\bibitem{Be77}
E. Bertini, \emph{Ricerche sulle transformazioni univoche involutori del piano}, Ann. Mat. Ser. II \textbf{8} (1877), 224--286.%

\bibitem{BoTsch99}
F.~Bogomolov, Yu.~Tschinkel, \emph{On the density of rational points on elliptic fibrations}, J. Reine Angew. Math. \textbf{511} (1999), 87--93.%

\bibitem{BoTsch00}
F.~Bogomolov, Yu.~Tschinkel, \emph{Density of rational points on elliptic $K3$ surfaces}, Asian J. Math. \textbf{4} (2000), 351--368.%

\bibitem{CalLy94}
F.~Call, G.~Lyubeznik, \emph{A simple proof of Grothendieck's theorem on the parafactoriality of local rings}, Contemp. Math. \textbf{159} (1994), 15--18.%

\bibitem{CaCe82}
F.~Catanese, G.~Ceresa, \emph{Constructing sextic surfaces with a given number $d$ of nodes}, J. Pure Appl. Algebra \textbf{23} (1982), 1--12.%

\bibitem{Ch00a}
I.~Cheltsov, \emph{Log pairs on birationally rigid varieties}, Jour. Math. Sciences \textbf{102} (2000), 3843--3875.%

\bibitem{Ch00b}
I.~Cheltsov, \emph{On smooth quintic 4-fold}, Mat. Sbornik \textbf{191} (2000), 139--162.%

\bibitem{Ch00c}
I.~Cheltsov, \emph{Log pairs on hypersurfaces of degree $N$ in $\mathbb{P}^{N}$}, Math. Notes \textbf{68} (2000), 113--119.%

\bibitem{Ch01a}
I.~Cheltsov, \emph{A Fano 3-fold with unique elliptic structure}, Mat. Sbornik \textbf{192} (2001), 145--156.%

\bibitem{Ch03b}
I.~Cheltsov, \emph{Non-rationality of a four-dimensional smooth complete intersection of a quadric and a quartic not containing a plane}, Mat. Sbornik \textbf{194} (2003), 95--116.%

\bibitem{Ch04a}
I.~Cheltsov, \emph{Conic bundles with big discriminant loci}, Izviestia RAN (2004), toappear.%

\bibitem{ChPa04}
I.~Cheltsov, J.~Park, \emph{Sextic double solids}, arXiv:math.AG/0404452 (2004).%

\bibitem{ChWo04}
I.~Cheltsov, L.~Wotzlaw, \emph{Non-rational complete intersections}, Proc. Steklov Inst. (2004), to appear.%

\bibitem{ClGr72}
H.~Clemens, P.~Griffiths, \emph{The intermediate Jacobian of the cubic threefold}, Ann. of Math. \textbf{95} (1972), 73--100.%

\bibitem{Co95}
A.~Corti, \emph{Factorizing birational maps of threefolds after Sarkisov}, J. Alg. Geometry \textbf{4} (1995), 223--254.%

\bibitem{Co00}
A.~Corti, \emph{Singularities of linear systems and 3-fold birational geometry}, L.M.S. Lecture Note Series \textbf{281} (2000), 259--312.%

\bibitem{CPR}
A.~Corti, A.~Pukhlikov, M.~Reid, \emph{Fano 3-fold hypersurfaces}, L.M.S. Lecture Note Series \textbf{281} (2000), 175--258.%

\bibitem{CosDo89}
F.~Cossec, I.~Dolgachev, \emph{Enriques Surfaces I} (Birkhauser, Boston 1989).%

\bibitem{Do66}
I.~Dolgachev, \emph{Rational surfaces with a pencil of elliptic curves}, Izv. Akad. Nauk SSSR Ser. Mat. \textbf{30} (1966), 1073--1100.%

\bibitem{Do82}
I.~Dolgachev, \emph{Weighted projective varieties}, Lecture Notes in Math. \textbf{956} (1982), 34--71.%

\bibitem{En99}
S.~Endra\ss, \emph{On the divisor class group of double solids}, Manuscripta Math. \textbf{99} (1999), 341--358.%

\bibitem{FuLa81}
W.~Fulton, R.~Lazarsfeld, \emph{Connectivity and its applications in algebraic geometry}, Lecture Notes in Math. \textbf{862} (1981), 26--92.%

\bibitem{Gro65}
A.~Grothendieck et al., \emph{Cohomologie locale des faisceaux coh\'erents et th\'eor\'emes de Lefschetz locaux et globaux}, S\'eminaire de g\'eom\'etrie alg\'ebrique, IHES (1965).%

\bibitem{Hal82}
G.~Halphen, \emph{Sur les courbes planes du sixieme degre а neuf points doubles}, Bull. Soc. Math. France \textbf{10} (1882), 162--172.%

\bibitem{HaTsch00}
J.~Harris, Yu.~Tschinkel, \emph{Rational points on quartics}, Duke Math. J. \textbf{104} (2000), 477--500.%

\bibitem{Har80}
R.~Hartshorne, \emph{Stable reflexive sheaves}, Math. Ann. \textbf{254} (1980), 121--176.%

\bibitem{IF00}
A.~R.~Iano-Fletcher, \emph{Working with weighted complete intersections}, L.M.S. Lecture Note Series \textbf{281} (2000), 101--173.%

\bibitem{Is80}
V.~Iskovskikh, \emph{Birational automorphisms of three-dimensional algebraic varieties}, J. Soviet Math. \textbf{13} (1980), 815--868.%

\bibitem{Is97}
V.~Iskovskikh, \emph{On the rationality problem for three-dimensional algebraic varieties}, Proc. Steklov Inst. Math. \textbf{218} (1997), 186--227.%

\bibitem{IsMa71}
V.~Iskovskikh, Yu.~Manin, \emph{Three-dimensional quartics and counterexamples to the L\"uroth problem}, Mat. Sbornik \textbf{86} (1971), 140--166.%

\bibitem{Ish82}
Sh.~Ishii, \emph{A characterization of hyperplane cuts of a smooth complete intersection}, Proc. Japan Acad. \textbf{58} (1982), 309--311.%

\bibitem{JaRu97}
D.~Jaffe, D.~Ruberman, \emph{A sextic surface cannot have $66$ nodes}, J. Alg. Geometry \textbf{6} (1997), 151--168.%

\bibitem{Kaw01}
M.~Kawakita, \emph{Divisorial contractions in dimension $3$ which contracts divisors to smooth points}, Invent. Math. \textbf{145} (2001), 105--119.%

\bibitem{KMM}
Y.~Kawamata, K.~Matsuda, K.~Matsuki,  \emph{Introduction to the minimal model problem, }Adv. Stud. Pure Math. \textbf{10} (1987), 283--360.%




\bibitem{Ko91}
J.~Koll\'ar et al., \emph{Flips and abundance for algebraic threefolds }Ast\'erisque \textbf{211} (1992). %

\bibitem{Ko95}
J.~Koll\'ar, Journal of the AMS \textbf{8} (1995), 241--249.%

\bibitem{Ko96}
J.~Koll\'ar, \emph{Rational curves on algebraic varieties} (Springer, Berlin 1996).%

\bibitem{Ko97}
J.~Koll\'ar, \emph{Singularities of pairs}, Algebraic geometry---Santa Cruz 1995 AMS (1997), 221--287.%

\bibitem{Ko00}
J.~Koll\'ar, \emph{Non-rational covers of $\mathbb{CP}^{n}\times \mathbb{CP}^{m}$}, L.M.S. Lecture Note Series \textbf{281} (2000), 51--71.%

\bibitem{Ne71}
M.~Noether, \emph{Uber Fl\"achen, welche Shaaren rationaler Curven besitzen}, Math. Ann. \textbf{3} (1871), 161--227.%

\bibitem{Pet98}
K.~F.~Petterson, \emph{On nodal determinantal quartic hypersurfaces in $\mathbb{P}^4$}, Thesis, University of Oslo (1998).%

\bibitem{Pr01}
Yu.~Prokhorov, \emph{Lectures on complements on log surfaces}, MSJ Memoirs \textbf{10} (2001).%

\bibitem{Pu88a}
A.~Pukhlikov, \emph{Birational automorphisms of a double space and a double quartic}, Izv. Akad. Nauk SSSR Ser. Mat. \textbf{52} (1988), 229--239.%

\bibitem{Pu95}
A.~Pukhlikov, \emph{A remark on the theorem on a three-dimensional quartic of V.A.Iskovskikh and Yu.I.Manin about 3-fold quartic}, Proc. Steklov Inst. Math. \textbf{208} (1995), 244--254.%

\bibitem{Pu97}
A.~Pukhlikov, \emph{Birational automorphisms of double spaces with sigularities}, J. Math. Sci. \textbf{85} (1997), 2128--2141.%

\bibitem{Pu00}
A.~Pukhlikov, \emph{Essentials of the method of maximal singularities}, L.M.S. Lecture Note Series \textbf{281} (2000), 73--100.%

\bibitem{Re00}
M.~Reid, \emph{Graded rings and birational geometry}, Proc. Symp. Alg. Geom. (Kinosaki) (2000) 1--72.%

\bibitem{Ry02}
D.~Ryder, \emph{Elliptic and K3 fibrations birational to Fano 3-fold weighted hypersurfaces}, Thesis, University of Warwick (2002).%

\bibitem{Sh92}
V.~Shokurov, \emph{Three-dimensional log perestroikas}, Izv. Ross. Akad. Nauk Ser. Mat. \textbf{56} (1992), 105--203.%

\bibitem{St78}
E.~Stagnaro, \emph{Sul massimo numero di punti doppi isolati di una superficie algebrica di $\mathbb{P}^{3}$}, Rend. Sem. Mat. Univ. Padova \textbf{59} (1978), 179--198.%

\bibitem{vSt93}
D.~van~Straten, \emph{A quintic hypersurface in $\mathbb{P}^{4}$ with $130$ nodes}, Topology \textbf{32} (1993), 857--864.%

\bibitem{Va83}
A.~Varchenko, \emph{On semicontinuity of the spectrum and an upper bound for the number of singular points of projective hypersurfaces}, Dokl. Aka. Nauk USSR \textbf{270} (1983), 1294--1297.%

\bibitem{Wa98}
J.~Wahl, \emph{Nodes on sextic hypersurfaces in $\mathbb{P}^3$}, J. Diff. Geom. \textbf{48} (1998), 439--444.%

\bibitem{Zak93}
F.~Zak, \emph{Tangents and secants of algebraic varieties}, Mathematical Monographs \textbf{127} (AMS, Providence 1993). %

\bibitem{Za52}
O.~Zariski, \emph{Complete linear systems on normal varieties and a generalization of a lemma of Enriques--Severi}, Ann. of Math. \textbf{55} (1952), 552--552.%

\bibitem{Za58}
O.~Zariski, \emph{On Castelnuovo's criterion of rationality $p_{a}=P_{2}=0$ of an algebraic surface}, Illinois J. Math. \textbf{2} (1958), 303--315.%

\end{thebibliography}
\end{document}